\documentclass[12pt]{amsart}

\usepackage[a4paper]{geometry}

\usepackage{amsfonts, amsmath, amssymb, amscd, latexsym, graphicx, pb-diagram, caption, mathdots}

\def\nn{\mathbb{N}}

\def\rr{\mathbb{R}}

\def\qed{\hfill\rlap{$\sqcup$}$\sqcap$\par}

\newtheorem{thm}{Theorem}[section]
\newtheorem{prop}[thm]{Proposition}

\newtheorem{defn}[thm]{Definition}

\newtheorem{question*}{Question}

 \theoremstyle{definition}
  \newtheorem{example}[thm]{Example}

\title{Grid diagrams for higher-dimensional Thompson's groups}
\author{Jos\'e Burillo }
\address{Departament de Matem\`atiques,
Universitat Polit\`ecnica de Catalunya, Jordi Girona 1-3, 08034
Barcelona, Spain} \email{pep.burillo@upc.edu}
\thanks{The first author thanks the Spanish Ministry MICINN through grant PID2021-126851NB-I00P for their support. }
\author{Sean Cleary}
\address{The City College of New York \and The CUNY Graduate Center, New York, NY 10031, USA}
\email{scleary@ccny.cuny.edu} \thanks{Support for this project was provided by a PSC-CUNY Award, jointly funded by The Professional Staff Congress and The City University of New York.}

\author{Brita Nucinkis}
\address{Department of Mathematics, Royal Holloway, University of London,
Egham, TW20 0EX, UK}
\email{Brita.Nucinkis@rhul.ac.uk}

\date{}

\begin{document}

\maketitle

\begin{abstract}
We describe standard forms for elements of the higher-dimensional Thompson groups $nV$ arising from gridding subdivision processes.
These processes lead to standard normal form descriptions for elements in these groups, and sizes of these standard forms estimate the word length with respect to
finite generating sets.  These gridded forms lead to standard algebraic descriptions as well, with respect to the both infinite and finite generating sets for these groups.

\end{abstract}

\section*{Introduction }

Thompson's groups form a remarkable family and arise in a range of settings.  Recurring themes within the family are the tension between working with infinite generating sets, which are convenient for representing elements, and finite generating sets, where the metric properties are clearer.  A common approach
for groups in the Thompson family for representing elements with respect to infinite generating sets is to factor elements into appropriate pieces, such as positive and negative parts, which give useful normal forms with both geometric and algebraic understanding.
  In the case of  Thompson's group $F$, there are tight correspondences between the geometric and algebraic normal forms, which give effective ways of understanding some of the properties of $F$.  In the cases of $T$ and $V$, there are normal forms and standard geometric representations, which are useful for understanding group properties.    Other groups in the Thompson family such as the braided Thompson's groups have normal forms that have geometric and algebraic descriptions, but the connections are not as strong as in the case of $F$.

 The higher-dimensional Thompson's groups $nV$, introduced by Brin \cite{brinnV} have descriptions which are more graphical in nature, typically described as bijections between dyadic dissections of $n$-dimensional cubes.  Here, we describe some normal forms for groups in this family as having canonical normal forms arising from a subdivision process we term ``gridding'' which give a means of understanding group properties, metric properties, and which leads to algebraic normal forms in these groups.

\section{Background}

Normal forms for Thompson's groups play an important role in describing and understanding group elements.  The normal forms for Thompson's group $F$ with
respect to the infinite generating set, described by Brown and Geoghegan \cite{bg} are exceptionally useful and form an ideal model
for normal forms for other members of the Thompson family.  These normal forms are of the type  $$x_{i_1}^{r_1}
x_{i_2}^{r_2}\ldots x_{i_k}^{r_k} x_{j_l}^{-s_l} \ldots
x_{j_2}^{-s_2} x_{j_1}^{-s_1} $$
where all of the exponents satisfy $r_i, s_i >0$, and the indices are arranged so that $i_1<i_2
\ldots < i_k$ and $j_1<j_2 \ldots < j_l$, with the additional condition that if  both $x_i$ and $x_i^{-1}$
occur, so does at least one of $x_{i+1}$ or  $x_{i+1}^{-1}$.  This standard normal form omits generators with zero exponent for brevity, but there is also an alternative approach to normal forms which include all relevant generators in appropriate range, including those occurring with zero exponent.  These
normal forms are $$x_{0}^{r_0}
x_{1}^{r_1}\ldots x_{k}^{r_k} x_{l}^{-s_l} \ldots
x_{1}^{-s_1} x_{0}^{-s_0} $$
where all of the exponents satisfy $r_i, s_i  \geq 0$, and the additional conditions are now that if  both $r_i$ and $s_i$ are non-zero, at least one of $r_{i+1}$ or $s_{i+1}$ must be non-zero, as well as requiring that the exponents  of the highest subindex generators (that is, $r_k$ and $s_l$)  are nonzero. These normal forms include the zero-exponent generators and thus can be longer but have useful related properties in the algebraic descriptions of normal forms in $nV$ below.

Either of these normal forms can be described as having a form $P N^{-1}$ where $P$ and $N$ are positive words, and have algebraic and geometric meaning and have a strong connection via leaf exponents, described by Cannon, Floyd, and Parry \cite{cfp}.  The algebraic forms give an effective way to multiply elements using the defining relations, and the geometric forms raise the recurring theme of factoring group elements into positive and negative factors.  That is, we can regard the positive part of the normal form as described by a tree pair diagram where the target tree is an all-right tree, and the negative part is described as a tree pair diagram where the source tree is an all-right tree.  There are corresponding equally effective normal forms for the groups $F(p)$, see Burillo, Cleary, and Stein \cite{bcs}.

For Thompson's groups $T$ and $V$ there are again algebraic and geometric normal forms for elements.  As described in Cannon, Floyd and Parry \cite{cfp}, in $T$, the normal forms are of the type for $F$ with an additional  middle term corresponding to a power of cyclic rotation on an all-right tree of an appropriate size, giving
 $$x_{i_1}^{r_1}
x_{i_2}^{r_2}\ldots x_{i_k}^{r_k} C_l^k x_{j_l}^{-s_l} \ldots
x_{j_2}^{-s_2} x_{j_1}^{-s_1} $$ which we regard as s $P C N^{-1}$ where $P$ and $N$ are as above, and $C$ is a power of one of the cyclic generators $c_n$.
Geometrically, these forms correspond to factorizations of a group element as a positive tree pair from an arbitrary tree to an all-right tree, then a cyclic term which reindexes the  leaves of an all-right tree by a cyclic shift of leaf indices, and then a negative tree pair from an all-right tree to an arbitrary tree.
 In the case of $V$, we have similar normal forms  $P \Pi N^{-1}$, where $\Pi$ is a permutation of the leaves of an all-right tree.  Note that although this gives some normal forms, the algebraic form is not nearly as canonical as the cases for $F$ and $T$ as there may be many ways to represent the permutation $\Pi$ in terms of the conjugates of transpositions of leaves of the all-right tree.  For braided Thompson's groups, there are natural normal forms of the form $P Br N^{-1}$, where $Br$ is a braid on the strands of the all-right tree, see the works of Brin \cite{brinBV1,brinBV2} and Dehornoy \cite{dehornoy}, and similarly for the pure braided Thompson's group F, normal forms  $P (Pr) N^{-1}$, where $Pr$ is a pure braid on the leaves of the all-right tree, see Brady, Burillo, Cleary, and Stein \cite{bbcs}. %We note that these normal forms often do not give effective insight into the metric properties, ...

 The descriptions of the Thompson's groups so far all have a natural left-to-right order on the leaves corresponding to the intervals in the dyadic subdivisions.
 These left-to-right orders lead themselves to natural descriptions by trees.  The higher dimensional versions $nV$ of Thompson's group $V$  do not have such
 a natural left-to-right order on the dyadic
subdivisions arising in the partitions, but nevertheless have
effective descriptions in terms of trees, as
 described by
 Burillo and Cleary \cite{2vm}, and Kochloukova,  Mart{\'{\i}}nez-P{\'e}rez,  and Nucinkis \cite{cohom2V3V}.  For example, in $2V$, the carets of a tree indicating a natural ``halving" of intervals to make the subdivisions  have corresponding vertical or horizontal ``halving" procedures, where the carets are now of different possible shapes (triangular, square, rounded, etc.), colors, dashings, or marked with labels to indicate in which direction of the multiple possible dimensions the halving takes place.   Below, we will use this to give geometric normal forms for elements of $nV$ of a ``multi-caret-type tree / permutation /  multi-caret-type tree'' type, and to those we can associate  algebraic normal forms as well.

 First, we generalize the notion of a basic dyadic interval to that of a basic dyadic block.
\begin{defn} A subinterval of the unit interval is called a \emph{basic dyadic
interval} if there exists $k\in\nn$, and also
$i\in\{1,2,\ldots,2^k\}$ such that the interval is
$$
\left[\frac{i-1}{2^k},\frac{i}{2^k}\right].
$$
\end{defn}
Basic dyadic intervals correspond to the vertices of the infinite
binary tree obtained by iterated subdivisions the unit interval in
halves.
\begin{defn} We will call a \emph{basic dyadic block} a subset of $\rr^n$ which is
the product of $n$ basic dyadic intervals. In particular, a
2-dimensional basic dyadic block will be called a \emph{basic dyadic
rectangle.}
\end{defn}

\section{Grid tree-pair diagrams in $2V$}

Though these forms apply to all higher dimensional groups $nV$, we
first study $2V$ for simplicity, and introduce  \emph{grid
tree-pair diagrams}. We start with some definitions.

\begin{defn} A \emph{grid subdivision} on a unit square is a subdivision
of the unit square such that each subdividing segment runs the full
length of the square (spanning several rectangles if necessary). A unit
square with a grid subdivision is a \emph{grid square.}
\end{defn}

\begin{figure}[h]
\centerline{\includegraphics[width=40mm]{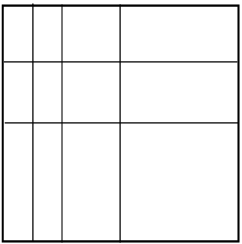}} \caption{A grid
square, which can be regarded as a product subdivision of two dyadic subdivisions of the horizontal and vertical directions.} \label{square}
\end{figure}

That is, a grid subdivision has the property that if a dyadic coordinate is the starting or ending coordinate of some rectangular block with a value in the other coordinate, that
starting or ending coordinate also appears in with all other possible other values in the other coordinate.   Thus, we can regard these grid subdivisions as products of dyadic subdivisions of the
two coordinates.
The square is subdivided into a grid of rectangles, hence the
name. An example is illustrated in Figure \ref{square}.

\begin{defn} A \emph{grid tree-pair diagram} (or just \emph{grid
diagram} for short) is a diagram where the
source square is a grid square.
\end{defn}

We begin by showing that every element has a grid diagram as a representative.

\begin{prop} Each element in $2V$ admits a grid tree-pair diagram.  \label{hasagrid}
\end{prop}

{\it Proof.}
To obtain a grid tree-pair diagram, we take all segments appearing in the subdivision of
the square and extend them until they run the full length of the
square. Obviously this means some other rectangles are subdivided,
so we carry these subdivisions forward via the element to perform the corresponding subdivisions in the target square, giving a grid-tree pair representative
of the desired element.  \qed

It is in general unreasonable to expect that an element could have a diagram with
both source and target to be grid squares.  For example, once the
source diagram is gridded, the process of gridding the resulting target arrangement
of rectangles will almost always result in subdivisions in the source which do not
cross the entire grid, undoing the gridding of the source subdivision.
There are a number of possible choices here, but for simplicity we
choose to take the source subdivision to be in gridded form.

\begin{figure}[ht]
\centerline{\includegraphics[width=120mm]{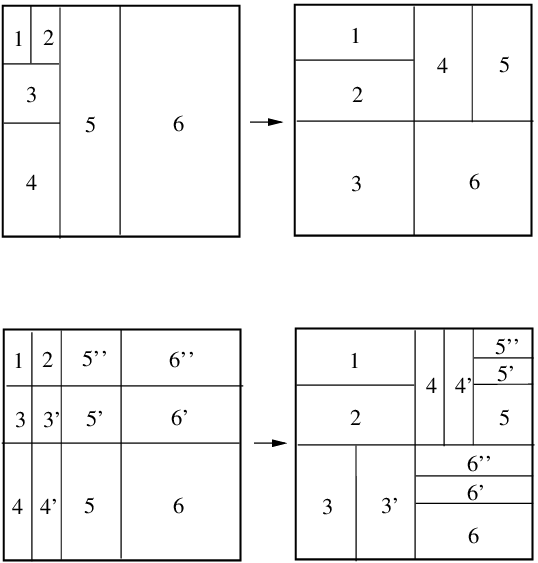}} \caption{Two representations of the same element $w$ of $2V$.  The top
representative is not gridded, but is refined to a grid tree-pair diagram in the bottom representation, where the source square is a product of
two dyadic subdivisions of the interval.}
\label{griddiagram}
\end{figure}

We illustrate in Figure \ref{griddiagram}  an example of an element $w$ and its resulting
grid diagram from the process described above.

\begin{figure}[p]
\centerline{\includegraphics[width=90mm]{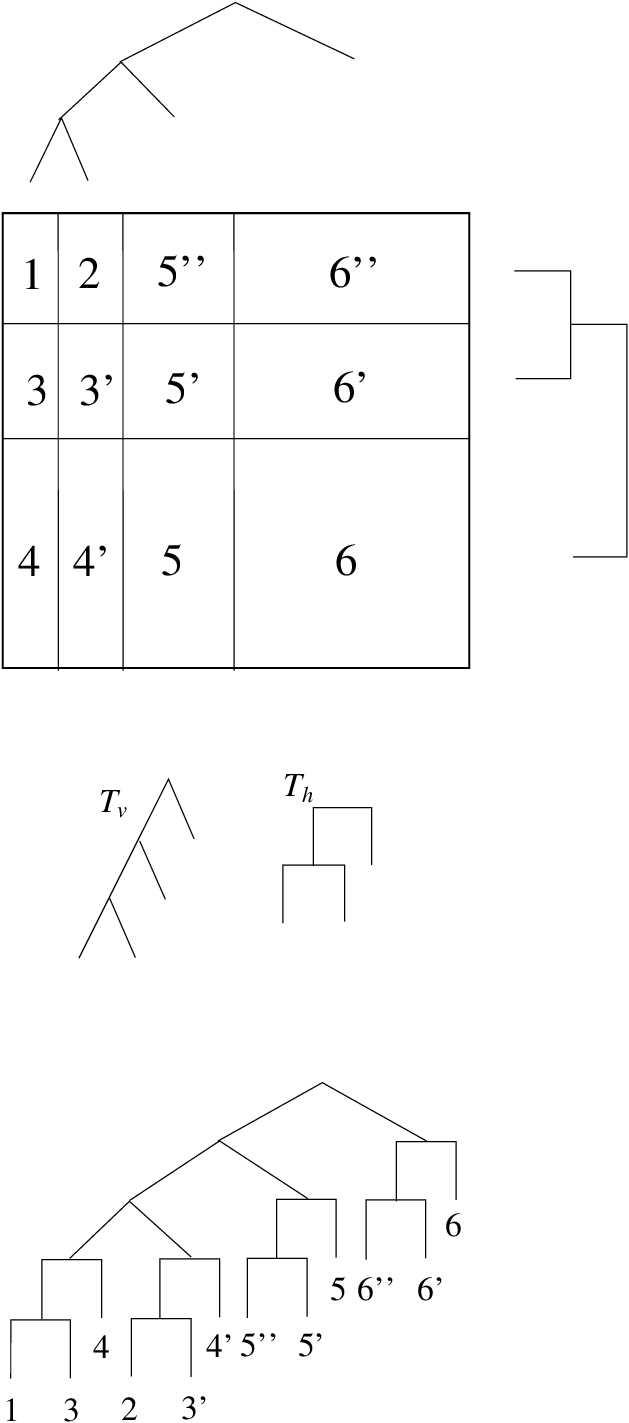}} \caption{The
source square for the example element $w$ and the two trees that
define it. The tree $T_h$ is represented on the right of the square
because left leaves represent top halves in subdivisions. At the
bottom, a tree representing this subdivision where copies of $T_h$ hang from
the leaves of $T_v$. An analogous tree representative for the same
grid subdivision could be constructed with copies of
$T_v$ hanging from the leaves of $T_h$, with the same dyadic blocks occuring, albeit in a different order.} \label{gridsquare}
\end{figure}

Grid squares can be understood also in terms of trees. Clearly, the
grid square can be codified with two trees: one tree representing the dyadic vertical
subdivisions and having only triangular, vertical-type, carets, and another tree representing the dyadic horizontal subdivisions, which will have only square, horizontal-type,
carets, using the convention described in \cite{2vm}.  We illustrate this process in  Figure \ref{gridsquare}.

So clearly two trees like these define a unique grid square and
vice versa. The actual tree for the square has one copy of the two trees
at the root, and many copies of the other tree, each copy hanging from a leaf of the first tree. This
can be done with the vertical tree at the root and the horizontal at
the leaves, or also the other way around.  The left-to-right order of leaves in the two orders will differ,
as commonly seen in equivalent tree-based representatives of elements of higher-dimensional Thompson's groups, but
the sets of dyadic blocks will coincide and the potential differing orders will result in a different permutation on the leaves to represent
the same group element.

As an example, the standard relation in $2V$ expressing the simplest interaction between iterated horizontal
and vertical subdivisions can be represented with
a grid square with one vertical and one horizontal (full-length)
subdivision, so both the horizontal and vertical trees used to build the tree will each have a single caret.  We assemble these trees in the
two possible manners, first hanging two copies of the vertical trees (each a single caret) from each of the two leaves of a single horizontal caret,
and then hanging two copies of the horizontal trees (each a single caret) from each of the two leaves of a vertical caret, giving the relation pictured in Figure 4 in \cite{2vm} with a permutation to match the corresponding dyadic blocks.

\begin{figure}[p]
\centerline{\includegraphics[width=130mm]{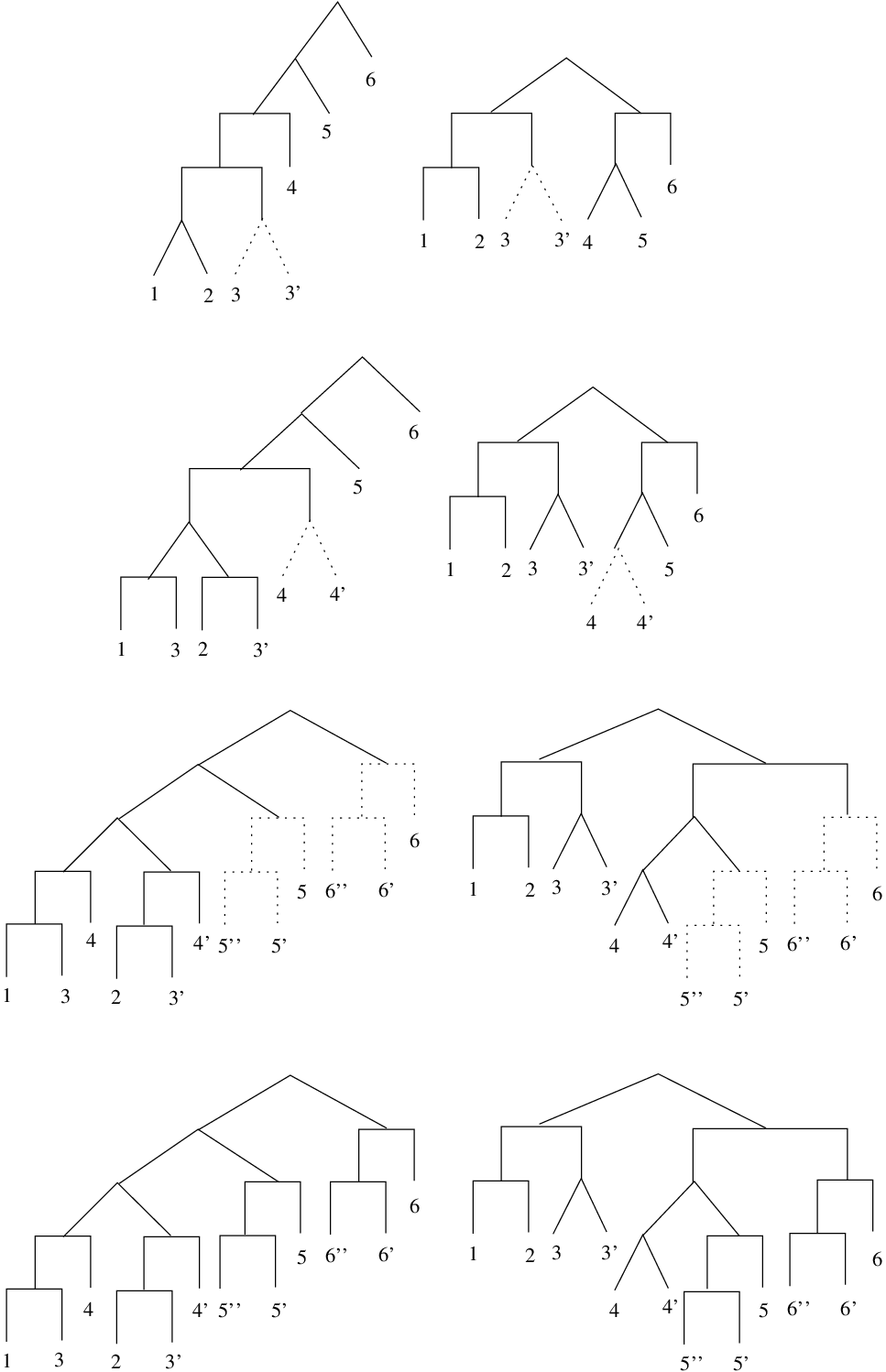}} \caption{The
process of transforming a diagram of our element $w$ into a grid
tree-pair diagram. The top tree-pair corresponds to the element $w$
depicted in Figure \ref{griddiagram}. The process consists of adding
vertical-type carets (of triangular type, in dotted lines) needed to apply the defining
relation for $2V$ and bring vertical carets to the top to form the
tree $T_v$. Finally, some horizontal carets (of square type, again dotted) have been added to have full copies of
$T_h$ hanging from each leaf of $T_v$.} \label{manytrees}
\end{figure}

The process of taking a subdivided square and finding a grid square for
it can also be seen with trees. From a starting tree representative, suppose
the desire is  to have all vertical carets on top and horizontal ones below those.
So we may add vertical carets to all leaves that need
them so that vertical carets appearing in the tree can be moved up.
We illustrate in Figure \ref{manytrees} how this process applies to the example element $w$ from Figures \ref{gridsquare} and \ref{manytrees}.

\section{Reduced grid diagrams}

We want to use grid diagrams to extend the notion of \emph{reduced
diagrams} we have in standard Thompson's groups (such as $F$, $T$ or
$V$) to $2V$. Recall that with relation to standard diagrams such as
those in \cite{2vm}, elements of $2V$ can have more than one reduced
tree-pair diagram (reduced in the traditional sense). We will define
reduced grid diagrams and prove their uniqueness.

We let $f$ be an element of $2V$ and we consider a grid diagram for
it. The source square is a grid square, so there are two trees $T_v$
and $T_h$ which represent the grid square via the product process
described above. The tree  $T_v$ has only vertical-type carets and
$T_h$ has only horizontal-type carets. We consider the basic dyadic
intervals $I_1,I_2,\ldots,I_m$ and $J_1,J_2,\ldots,J_n$ which
represent the ordered leaves of the two trees; namely, the $I_i$ are
intervals on the bottom edge of the square and the $J_j$ are on the left
edge. So each rectangle of the grid is of the type $I_i\times J_j$,
and it naturally is a basic dyadic rectangle.

We observe how there could be a possible reduction of the tree
$T_v$. We would be able to erase a caret of it if its two leaves are
represented by intervals $I_i$ and $I_{i+1}$ which could be joined
into a single interval. This could only be done if the map $f$ is
completely affine in all the rectangles $(I_i\cup I_{i+1})\times
J_j$ for every $j=1,2,\ldots,n$, and the image of each of these
basic rectangles is also a basic dyadic rectangle. So we can say
that the element is {\em vertically reduced} if no caret in $T_v$
can be reduced. That is, for every exposed caret of $T_v$ and its
leaves $I_i$ and $I_{i+1}$, there exists a $j\in\{1,2,\ldots,n\}$
such that the affine map is different on $I_i\times J_j$ and in
$I_{i+1}\times J_j$. Or perhaps it is the same affine map, but the
image is not a basic dyadic rectangle.

This notion can be understood easily: a diagram is vertically
nonreduced if for an exposed caret in $T_v$, the whole segment in
the square that separates its two leaves could be erased and the map
would be the same (and the image of a basic dyadic rectangle is a
basic dyadic rectangle), and naturally we say a diagram is {\em
vertically reduced} if there are no possible vertical reductions.

Clearly we can define an element to be {\em horizontally reduced}
analogously  for the carets on $T_h$ and its rectangles. This leads
to the definition of reduced grid tree-pair diagrams.

\begin{defn} A grid tree-pair diagram for an element of $2V$ is
called \emph{reduced} if it is both vertically and horizontally
reduced.
\end{defn}

This just means that the two trees $T_v$ and $T_h$ cannot be
reduced, so they are the smallest possible trees which can be used to
describe the element via this product subdivision process. If a caret is removed in
one of them, then the map $f$ cannot be represented with the grid
square corresponding to these two trees.

The main theorem is the following:

\begin{thm} An element $f$ in $2V$ has a \emph{unique} reduced grid
tree-pair diagram. \label{uniqgrid}
\end{thm}

{\it Proof.} The reader is encouraged to look up the proof for the
uniqueness of reduced diagrams in $F$ found in Cannon, Floyd, and Parry \cite{cfp}, since our
proof follows the same lines.

We consider a grid diagram for $f$, and reduce the two trees as much
as possible. We let $I_1,I_2,\ldots,I_m$ and $J_1,J_2,\ldots,J_n$ be
the vertical and horizontal intervals in which the unit interval is
subdivided, vertically and horizontally. We observe that all $I_i$
and $J_j$ are basic dyadic intervals of length $1/2^k$ for some $k$.

So to understand uniqueness we only need to see that a basic dyadic
interval $I$ in the horizontal edge of the square has all the
intervals $I\times J_j$ mapped linearly by $f$ ($j=1,\ldots,n$) into
basic dyadic rectangles if and only if $I$ is represented in $T_v$
by a leaf, or else is not at all in $T_v$. If $I$ were represented
by an internal node in $T_v$, then at least one of the rectangles
$I\times J_j$ would be subdivided vertically into two rectangles
with different affine maps. Otherwise, the tree could be reduced
until $I$ was represented by a leaf.

This characterization of the tree $T_v$ requires the use of the
intervals $J_j$, which may appear to be a circular argument, since
$T_h$ would use the intervals $I_i$. But we note that this can be
done in a precise way using only the $I_i$. We can characterize the
intervals $I_i$ (without using the $J_j$) the following way. A basic
dyadic interval $I$ on the bottom interval corresponds to a node
(not a leaf) of $T_v$ if and only if for every $n\in\nn$ there
exists a basic dyadic interval $K$ of length $1/2^n$ such that
$I\times K$ is not mapped affinely into a basic dyadic rectangle.
Equivalently, an interval $I$ in the bottom edge corresponds to a
leaf of $T_v$, or is not at all in $T_v$, if and only if there
exists $n_0\in \nn$ such that for all $n\ge n_0$ and all dyadic
intervals $K$ of length $1/2^n$ the rectangle $I\times K$ is mapped
affinely into a basic dyadic rectangle.

Clearly the same characterization works for the intervals defining
the tree $T_h$. Hence the reduced trees are actually completely
determined by the map, so they are unique.\qed

We have developed this construction in $2V$ for the sake of
simplicity, but we observe that grid diagrams generalize readily to
$nV$. A grid cube will have a subdivision into blocks using parts of
hyperplanes, and the grid cube can be codified as a product of $n$ trees, each
of which has a single type of carets (of the $n$ available types of
carets in $nV$). The process goes the exact same way and the
theorems hold with proofs adapted to the $n$-dimensional
case.

\section{Metric properties}

For the metric properties of $nV$, we use the word length with respect to
the finite generating sets.  The word length is estimated by Burillo and Cleary
 \cite{2vm},
in terms of the number of carets, in  a \emph{minimal-size tree-pair diagram.}
Such a diagram, although
it exists, may be difficult to compute, since finding it could be as
difficult as finding a minimal word for the element.  Furthermore, such minimal-size tree-pair diagrams
are not unique as there may be multiple representatives of a given element, each of minimal size.

Since now we have a unique well-defined diagram for each element,
which is easy to compute directly, here we consider the metric
properties via the number of carets of the reduced grid diagram,
which is unique and algorithmically constructible in an easy way.

So we let $x$ be an element of $2V$, and we denote by $|x|$ its word
length with respect to  the standard finite generating set for $2V$. As in \cite{2vm},
we let $N(x)$ be the number of carets in a \emph{minimal-size} reduced
tree diagram for $x$. And finally, we let $M(x)$ be the number of
carets in the reduced grid diagram for $x$. The result in \cite{2vm}
was that
$$
\log N(x) \prec |x| \prec N(x)\log N(x)
$$
where $\prec$ means the usual bounds  up to multiplicative and additive constants as arise in quasi-isometric embeddings.
So we note now that obviously $N(x)\leq M(x)$ (since $N$ is the
minimal number of carets), giving $M(x)$ as an obvious upper bound for $N(x)$.

For a lower bound, we note that $M(x)\leq N(x)^2$ by the
construction of the grid diagram. To see this, we take the minimal
diagram and observe that there are at most $N$ vertical subdivisions
and also at most $N$ horizontal ones. So the grid diagram can be
constructed with at most $N^2$ rectangles. Then, since $\log N(x)$
and $\log M(x)$ differ at most by a multiplicative constant 2, we
have that the lower bound still holds with adjusted linear coefficients.
 Hence, $M$ satisfies the same
inequality as $N$ for estimating word length via the number of
carets, giving the following theorem.

\begin{thm} The word length $|x|$ of an element $x\in 2V$, and the
number of carets $M(x)$ of its reduced grid diagram satisfy
$$
\log M(x) \prec |x| \prec M(x)\log M(x).
$$ \label{gridmetric}
\end{thm}

These bounds cannot be easily improved. The lower bound has the same
properties as it did for $N(x)$, since the elements representing the
exponential distortion are already represented by grid diagrams (see
\cite{2vm}). For the upper bound, one could speculate that a better
quantity to use could be the maximum of the number of carets of
$T_h$ and $T_v$, call this number $m(x)$. But this number would
roughly be $m(x)=\sqrt{M(x)}$, and it would not satisfy an analogous
upper bound: the number of elements with a given $M$ is at least
$M!$ due to the many possible permutations for those subdivisions, whereas the proposed upper bound
$m(x)\log m(x)$ would yield approximately
$$
\exp{(m\log m)}=\exp{(\sqrt M \log M)} = M^{\sqrt M}
$$
which is smaller than $M!$.  In fact, following arguments parallel to those of Birget \cite{birget}, we can see that the factorially many possible permutations
of the number of blocks in a subdivision means that most of the elements will have sizes quite close to the upper bound, as was shown by Birget in the case of $V$.

\section{Algebraic Normal Forms for Elements of $2V$ \label{algsec}}

We note that though the grid normal form is principally a geometric
normal form for elements of $nV$, there are associated unique algebraic
normal forms.   These arise in the following manner, described first
for $2V$ and then more broadly in the following section.  As in \cite{2vm}, for 2V, we
describe vertical divisions of the associated rectangles with
triangular carets, and horizontal divisions via square carets.

In what follows, we will use the standard infinite generating sets $\{A_i, B_i, C_i, \Pi_i, \overline{\Pi_i}\}$ described by Brin \cite{brinnV} for
$2V$. We use the $A_i$ generators for the triangular caret forms
of the standard Thompson's group generators where both the source
and the target are triangular, see also \cite{2vm}.  We use $B_i$ for the square caret generator forms,
where the source is of the square type and the target is of the
triangular type.   We use $C_i$ for the generators which have the
last ($i$th) caret of the all-right tree of square type, all other
carets in the source tree being triangular, and with the target tree
being an all-right tree with entirely triangular types.

We suppose that we have a normal form for the grid subdivision
of the source square
in the manner described in the previous section.  For the triangular subdivision tree $T_v$
associated to the vertical gridding, we consider a normal form $V$
associated to the tree that goes from that subdivision to the
all-right tree:

$$V=A_0^{a_0}A_1^{a_1}\ldots A_k^{a_k} $$
where all of the exponents satisfy $a_i\geq 0$, and $k+1$ is the
number of leaves for $T_v$. For the domain tree, the relevant leaf
indices are those from the positive part of the word, and in the
case where the gridding occurs entirely in the domain square, we
focus on that positive part of the algebraic form,
and will write $A_0^{a_0}A_1^{a_1}\ldots A_k^{a_k}$ for that. %where the $w$ is the negative generator part of the word and describes the target part of the subdivision.

Similarly, for the squared caret horizontal subdivision, we have a
normal form $H$ for the tree $T_h$ of the type $$H=B_0^{b_0}
B_1^{b_1}\ldots B_{k'}^{b_{k'}} $$ where again all of the exponents
satisfy $b_{i'}\geq 0$.

One complication is that the generators corresponding non-triangular carets on the spine
of the tree are created in a different manner than those on the
interior, so that affects the algebraic normal form construction.

We begin by describing what happens to the generators given by the
$T_h$-trees hanging off a leaf of $T_v$ that is not the right-most
one. Notice that the word $H$ reappears inserted multiple times in
$V$ at each leaf location, with its subindices shifted to account
for the increased starting location of the first leaf index.  We use
$H(k)$ to indicate the word obtained from $H$ by increasing each
subindex of the $y$ generators in $H$ by $k$, and we  see  a word of
the form $A_{i_1} H(...) A_{i_2} H(...) \ldots $ as a positive word
in the $A$- and $B$-generators for the gridded domain subdivision of
the square, where each of the $H(...)$ are the original word for $H$
with appropriate increases in the subindices for the $B$ generators
occurring there.

Note that the process of increasing subindices for generators in this manner
is a well-known process in $F$ and some of its generalizations that arises in many settings.  The shift homomorphism
$\phi: F \rightarrow F$ where $\phi(x_i)=x_{i+1}$ described by Brown and Geoghegan \cite{bg}
is an important homotopy idempotent (a universal example, in a sense). Here, the $x_i$ denote the generators in the standard infinite generating set for $F$.
Note that $H(k)= \phi^k(H)$ where $\phi$ is the shift automorphism
for the subgroup isomorphic to $F$ created by the generators $B_i$ in the horizontal direction.
 %We note that this important shift automorphism arises only in $F$ among the groups in the Thompson family, and not in $T$, $V$, $2V$ or necessarily not for any simple group in the family.

 What remains is dealing the with the generators given by $T_h$ hanging of the right-most leaf of $T_v$. This will involve the generators $B_i$ and $C_j.$ In particular, the smallest index $j$ will correspond to the number of carets on the spine of $T_v.$ The generators $C_j$ will, as in \cite{brinnV, 2vm}, appear on the left of the algebraic form. The $B_i$ will correspond to those carets not on the spine of $T_h,$ again shifted upwards depending on the number of leaves in the trees $T_v$ and $T_h$ as above. The word in the $B_j$ will then be behind the last $A$-generator.

This leads to a positive part $P$ of a normal form, built from the
gridded source diagram of the geometric normal form.  For the
remaining part of the normal form, there is the middle permutational
part $\Pi$ given by an appropriate product of conjugates of
transpositions of leaves of the all-right tree, followed by a
negative word  associated to a tree describing the target diagram of
the gridded diagram.  We note that are there few restrictions on the
permutational part or for the  negative word as there are many
possible elements with the same positive part of the normal form.
This gives an algebraic normal form  for the element of the type $P
\Pi N^{-1}$.

\begin{figure}[h]
\centerline{\includegraphics[width=60mm]{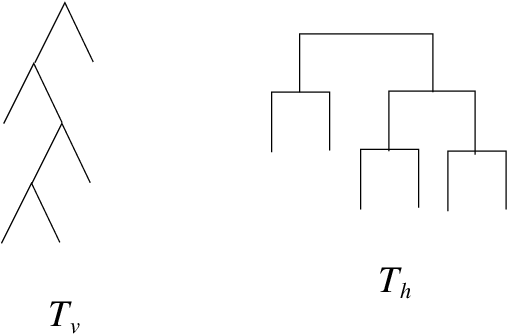}} \caption{Example
trees $T_v$ and and $T_h$ used to construct the positive part of the
algebraic normal form for an element $v$ as the vertical and
horizontal subdivisions in the product, respectively.} \label{tvth}
\end{figure}

\begin{example}  Here we illustrate this process
via an example element $f \in 2V$. We take a grid subdivision whose trees are the vertical and horizontal trees
following $T_v$ and $T_h$ respectively, illustrated in Figure \ref{tvth}.  As described above, the grid process
applies only to the source square, so we neglect in this example the permutational part and the domain subdivision
parts of $v$
as either of those could essentially be arbitrary, of the appropriate sizes and subject to not having obvious reductions present.

We construct the  the corresponding words for the positive part of  in the generators. For
$T_h$, this corresponds to the word $A_0A_1^2$. For convenience, we
are going to write this word in the form
$$
A_0^{a_0}A_1^{a_1}\ldots A_n^{a_n}
$$
where the exponent $a_i$ is allowed to be zero, and so we have
all generators listed, each corresponding to a leaf, even if that leaf does
not give a positive power for the relevant generator.
So our tree $T_v$ corresponds to
$A_0A_1^2A_2^0A_3^0A_4^0$. This will be convenient because this way
the generators are in bijection with the leaves of the tree.

We now  hang a copy of the tree $T_h$ at each leaf of
$T_v$ to obtain the source part of the grid diagram for $f$. This gives the tree illustrated in Figure \ref{bigtree}.

\begin{figure}[h]
\centerline{\includegraphics[width=140mm]{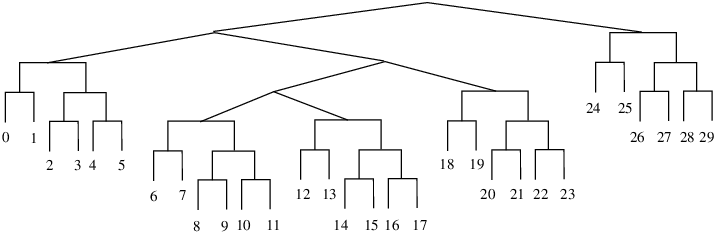}}
\caption{The tree used to construct the positive part of the normal form for $w$, obtained by attaching a copy of $T_h$ below each leaf of $T_v$.}
\label{bigtree}
\end{figure}

If a single copy of the tree $T_h$ were alone, the would be
represented as a positive element from the leaf exponents method by
$B_0B_1^0B_2B_3^0B_4B_5^0$, where again we choose to write zero
exponents for the generators not appearing. But when we hang this
tree from a leaf (except on the last leaf of $T_v$) all carets which
are on the spine of the tree now are in the interior, and hence they
account for positive powers of generators-- even the ones which did
not contribute earlier. Some exponents will be increased by 1. So
this first instance of the tree $T_h$ hanging from leaf 0 of $T_v$
contributes $B_0^2B_1^0B_2^2B_3^0B_4B_5^0$ to the algebraic form.

We inspect the algebraic word that is represented by this large tree of Figure \ref{bigtree}.  First, there is a power of $A_0$ followed by the word in the
$B$-generators given by $T_h$:
$$
A_0B_0^2B_1^0B_2^2B_3^0B_4B_5^0
$$
The continuation is a power of the next $A$-generator, which
corresponds to $A_1$ but here becomes $A_6$ since $T_h$ has 6 leaves and thus we have the increase in subindex of $A$.
Namely:
$$
A_0B_0^2B_1^0B_2^2B_3^0B_4B_5^0\ \
A_6^2B_6^2B_7^0B_8^2B_9^0B_{10}B_{11}^0
$$
This process continues with all the words corresponding to all the
copies of the tree $T_h$ except the final one. We obtain:
$$
\begin{array}c
A_0B_0^2B_1^0B_2^2B_3^0B_4B_5^0\ \
A_6^2B_6^2B_7^0B_8^2B_9^0B_{10}B_{11}^0\\
A_{12}^0B_{12}^2B_{13}^0B_{14}^2B_{15}^0B_{16}B_{17}^0\ \
A_{18}^0B_{18}^2B_{19}^0B_{20}^2B_{21}^0B_{22}B_{23}^0
\end{array}
$$
Finally, since the last copy of $T_h$ is located on the spine of the
big tree, it should be modelled on the word on the spine instead of
the interior word. Hence it should have some exponents diminished by
1, and also with the corresponding (exponent 0) $A$ at the
beginning. The resulting word for the positive part of $f$ is:
$$
\begin{array}c
A_0B_0^2B_1^0B_2^2B_3^0B_4B_5^0\ \
A_6^2B_6^2B_7^0B_8^2B_9^0B_{10}B_{11}^0\ \
A_{12}^0B_{12}^2B_{13}^0B_{14}^2B_{15}^0B_{16}B_{17}^0\\
A_{18}^0B_{18}^2B_{19}^0B_{20}^2B_{21}^0B_{22}B_{23}^0\ \
A_{24}^0B_{24}B_{25}^0B_{26}B_{27}^0B_{28}^0B_{29}^0
\end{array}
$$
This is not yet the word that corresponds to this tree, because the
horizontal carets on the spine (the $B$-generators that do not occur
 in the spine version of $T_h$) have to be created with
$C$-generators at the beginning of the word. So the exact word that
corresponds to this tree is
$$
\begin{array}c
C_1C_2C_3\ \ A_0B_0^2B_1^0B_2^2B_3^0B_4B_5^0\ \
A_6^2B_6^2B_7^0B_8^2B_9^0B_{10}B_{11}^0\ \
A_{12}^0B_{12}^2B_{13}^0B_{14}^2B_{15}^0B_{16}B_{17}^0\\
A_{18}^0B_{18}^2B_{19}^0B_{20}^2B_{21}^0B_{22}B_{23}^0\ \
A_{24}^0B_{24}B_{25}^0B_{26}B_{27}^0B_{28}^0B_{29}^0
\end{array}
$$

\end{example}

Note that this is just the positive part of the normal form, and the representative of the group element will be the positive part described above, followed by generators for the permutational part and then negative powers of the generators describing the tree for the target subdivision.

\begin{thm}\label{algform}

Every element $f \in V$ gives rise to a unique prefix word in the letters
$A_i, B_j, C_k$ arising from its reduced grid tree-pair diagram as
follows: Suppose the source grid gives rise to a vertical tree $T_v$
with $n+1$ leaves and a horizontal tree $T_h$ with $m+1$ leaves.
Then
\begin{itemize}
\item $T_v$  corresponds to a word
$A_0^{a_0}A_1^{a_1}\ldots A_n^{a_n}$, and
\item $T_h$ corresponds to two
words depending upon where in $T_v$ it hangs:
\begin{itemize}
\item $B_0^{b_0}B_1^{b_1}\ldots B_m^{b_m}$ when it is on the spine ($b_p$ can be zero),
\item $B_0^{b'_0}B_1^{b'_1}\ldots B_m^{b'_m}$ when it is {\bf not} on the
spine, (that is, if it hangs from the interior or left side of the
tree.) The exponent $b'_p$ is either equal to $b_p$ or $b_{p}+1$
depending on the generator.
\end{itemize}
\end{itemize}

Furthermore, the full tree corresponding to the source grid admits a unique word made out of the concatenation of
the following words:
\begin{itemize}
\item A word of the type $C_iC_{i+1}\ldots C_{i+\ell}$ corresponding to
the horizontal carets on the spine of $T_h$, where the indices are
all consecutive to create the right-hand side carets of the appropriate types,
\item $n$ words of the type
$$
A_{(m+1)j}^{a_j}B_{(m+1)j}^{b'_0}B_{(m+1)j+1}^{b'_1}\ldots
B_{(m+1)j+m}^{b'_m}
$$
for $j=0,1,\ldots,n-1$,
\item a final word, of the type
$$
A_{(m+1)n}^{a_n}B_{(m+1)n}^{b_0}B_{(m+1)n+1}^{b_1}\ldots
B_{(m+1)n+m}^{b_m}.
$$
\end{itemize}
\end{thm}

\begin{proof} This follows from  \cite[Theorem 2.2]{2vm}, Theorem \ref{uniqgrid},  and the considerations above.
\end{proof}

In the example above, we have $i=1$, $\ell=2$, $n=4$, and $m=5$.
The introduction of the generators with exponent zero helps
understand its regularity, but it makes the word unnecessarily long.
If we eliminate the terms with exponent 0, the word for the positive part of $v$ becomes just
$$
C_1C_2C_3\ \ A_0B_0^2B_2^2B_4\ \ A_6^2B_6^2B_8^2B_{10}\ \
B_{12}^2B_{14}^2B_{16}\ \ B_{18}^2B_{20}^2B_{22}\ \ B_{24}B_{26}
$$
which can be checked by reading exponents and generators directly from the tree.

Recall that these algebraic forms are with respect to the standard infinite generating sets $\{A_i, B_i, C_i, \Pi_i, \overline{\Pi_i}\}$ described by Brin \cite{brinnV}.  Using the standard expressions of the infinite generators in terms of the finite generating set $\Sigma$ described by Brin \cite{brinnV} as $\{A_0, A_1, B_0, B_1, \Pi_0, \overline{\Pi_0}, \Pi_1, \overline{\Pi_1}\}$ via substitution, the standard expressions with respect to the infinite
generating set give rise to standard expressions of elements with respect to the finite generating sets of the type $\Sigma$.  So grid tree diagrams give standard 
algebraic normal forms both with respect to infinite and finite generating sets.

\section{Other grid normal forms and higher dimensional cases}

In the higher dimensional cases of $nV$, there are similar grid diagrams arising from the partitions of the $n$-dimensional cube along each of its coordinates.  For the tree representative, we use $n$ different types of carets (possibly triangular, square, rounded, \ldots) and again can obtain a {\em grid subdivision}  which can be regarded as a product subdivision of  the appropriate dimensional unit cube via $n$ dyadic subdivisions of each of the $n$ coordinates.  There are natural generalizations of Proposition \ref{hasagrid} and Theorem \ref{uniqgrid} giving that such
elements have grid diagrams and that they are unique.    The metric estimates of Theorem \ref{gridmetric} also hold in $nV$ as instead of the number of carets growing quadratically at most, the number of carets in the gridding process grows at most as the $n$th power, but that potential growth still may be absorbed by adjusting the multiplicative constants outside of the  logarithms in the metric bound estimates.  The algebraic
form of Theorem \ref{algform} also give canonical normal forms, albeit more complicated as the natural ``nesting'' of the copies of the trees will lead to many terms.  For the positive part of the word, there will be the same components.  There will be a prefix word of generators of $n-1$ types to create the correct shape of the spine of the tree, followed by
a word in a standard form obtained by a nesting process.
Starting with words in the $A$ generators,  each gap between the powers of the $A$ generators will interrupted by subwords for the $B$ generators (with shifted indices, as above), and for which those $B$ subwords themselves will be interrupted by shifted index subwords of generators for each additional dimension in a natural recursive manner, nesting to give a standard normal form for the positive part of the word.

Also, we note that in $2V$, there was a choice to grid both of the coordinates in the source cube.  There are other possible choices, each of which
would lead to standard forms for that choice.   For example, we could have chosen to grid the target cube, or to grid the horizontal component in the source together with the vertical component in the target.
We can regard that decision as assigning to each coordinate a binary choice of that coordinate lying in either the source or target.  So there are 4 such possible choices for $2V$ and $2^n$ potential choices for $nV$, with each of those choices giving grid normal forms and associated algebraic normal forms.  Typically we presume that the normal forms are most useful when either the source or the target cube is completely gridded, rather than sharing the gridding between the two cubes.

If we consider subdivisions where the coordinates lie in
$\bold{Z}[\frac1p]$ for the corresponding versions of $nV$, these
methods give completely analogous results using multiple types of
carets with $p$ children in the place of the multiple types of
binary carets used above.  Similarly, for the groups $nBV$ where the
matchings between the subdivisions are given as braids in the
appropriate dimensional spaces rather than permutations, there are
similar standard forms for the positive and negative parts of the
word, with a middle word in a standard form chosen for the
appropriate braid representation.

\bibliographystyle{plain}

\def\cprime{$'$}

\end{document}